\documentclass[a4paper,12pt]{article}
\usepackage{times, url}
\usepackage{amsmath,amssymb,amsthm,amsfonts}

\theoremstyle{plain}
\newtheorem{thm}{Theorem}[section]

\newtheorem{lem}{Lemma}[section]

\theoremstyle{definition}

\newtheorem{dfn}{Definition}[section]
\theoremstyle{remark}

\numberwithin{equation}{section}

\begin{document}
\setcounter{page}{1}

\begin{center}
	
	{\LARGE \bf  Primality tests for specific classes of $N=k\cdot 2^m \pm 1$}
	\vspace{8mm}
	
	{\large \bf Predrag Terzi\' c}
	\vspace{3mm}
	
	Bulevar Pera \'Cetkovi\'ca 139 , Podgorica , Montenegro \\
	e-mail: \url{pedja.terzic@hotmail.com}
	\vspace{2mm}

\end{center}
\vspace{10mm}

\noindent
{\bf Abstract:}  Polynomial time primality  tests for specific classes of  numbers of the form $k\cdot 2^m \pm 1$ are introduced . 
\vspace{10mm}

\section{Introduction} 

In 1856 Edouard Lucas developed primality test for Mersenne numbers . The test was improved by Lucas in 1878 and Derrick Henry Lehmer in the 1930s , see \cite{1} .In 1960 Kusta Inkeri provided unconditional , deterministic , Lucasian type primality test for Fermat numbers , see \cite{2} .In 1969 Hans Riesel formulated primality test , see \cite{3}  for numbers of the form $k\cdot 2^n-1$ with $k$ odd and $k<2^n$ . In 2006 Zhi-Hong Sun provided general primality criterion for numbers of the form $k\cdot 2^m \pm 1$ ,see \cite{4} . In this note we present Lucasian type primality  tests for specific classes of $k\cdot 2^m \pm 1$ . 

Throughout this paper we use the following notations: $\mathbb{Z}$-the set of integers , $\mathbb{N}$-the set of positive integers , $\left(\frac{a}{p}\right)$-the Jacobi symbol , $(m,n)$-the greatest common divisor of $m$ and $n$ , $S_n(x)$-the sequence defined by $ S_0(x)=x $ and $ S_{k+1}(x)=(S_{k}(x))^2-2 (k\ge 0) $ .

\section{Basic Lemmas and Theorems}

\begin{dfn}\label{le}

For $P,Q \in \mathbb{Z}$ the Lucas sequence $\{V_n(P,Q)\}$ is defined by  \\
$V_0(P,Q)=2 ,V_1(P,Q)=P ,V_{n+1}(P,Q)=PV_n(P,Q)-QV_{n-1}(P,Q) (n\ge 1)$\\
Let $D=P^2-4Q$ . It is known that 
\begin{center}$V_n(P,Q)=\left(\frac{P+\sqrt{D}}{2}\right)^n+\left(\frac{P-\sqrt{D}}{2}\right)^n$
	\end{center}
	
	\end{dfn}
	
	\begin{lem}\label{le}
		
	Let $P,Q \in \mathbb{Z} $ and $n \in \mathbb{N} $. Then 
	\begin{center}
	$V_n(P,Q)=\displaystyle\sum_{r=0}^{[n/2]}\frac{n}{n-r}\binom{n-r}{r}P^{n-2r}(-Q)^r$	
		\end{center}
		
		\end{lem}

\begin{thm}\label{le}(Zhi-Hong Sun)
	
	For $m \in \{2,3,4,...\}$ let $p=k \cdot 2^m\pm 1$ with $0<k<2^m$ and $k$ odd . If $b,c \in \mathbb{Z}$ , $(p,c)=1$ and $\left(\frac{2c+b}{p}\right) = \left(\frac{2c-b}{p}\right)=-\left(\frac{c}{p}\right)$ then $p$ is prime if and only if $ p \mid S_{m-2}(x)$, where $x=c^{-k}V_k(b,c^2)=\displaystyle \sum_{r=0}^{(k-1)/2} \frac{k}{k-r}\binom{k-r}{r}(-1)^r(b/c)^{k-2r}$
	
\end{thm}

For proof see Theorem 3.1 in \cite{4} .

\begin{lem}\label{le}
	Let $n$ be odd positive number , then	
		\begin{center}
			$\left(\frac{-1}{n}\right) =
			\begin{cases}
				\hphantom{-}1, & \text{if } n \equiv 1 \pmod{4} \\
				-1, & \text{if } n \equiv 3 \pmod{4}
			\end{cases}	$
		\end{center}
		
\end{lem}

\begin{lem}\label{le}
	Let $n$ be odd positive number , then	
	\begin{center}
		$\left(\frac{2}{n}\right) =
		\begin{cases}
		\hphantom{-}1, & \text{if } n \equiv 1,7 \pmod{8} \\
		-1, & \text{if } n \equiv 3,5 \pmod{8}
		\end{cases}	$
	\end{center}
	
\end{lem}

\begin{lem}\label{le}
	Let $n$ be odd positive number , then\\
	case 1. ($n\equiv 1 \pmod 4$)	
	\begin{center}
		$\left(\frac{3}{n}\right)=
		\begin{cases}
		\hphantom{-}1&\text{if }n\equiv 1\pmod {3}\\
		\hphantom{-}0&\text{if }n\equiv 0\pmod {3}\\
		-1 & \text{if }n\equiv 2\pmod {3}
		\end{cases}$
	\end{center}
	case 2. ($n\equiv 3 \pmod 4$)
	\begin{center}
		$\left(\frac{3}{n}\right)=
		\begin{cases}
		\hphantom{-}1&\text{if }n\equiv 2\pmod {3}\\
		\hphantom{-}0&\text{if }n\equiv 0\pmod {3}\\
		-1 & \text{if }n\equiv 1\pmod {3}
		\end{cases}$
	\end{center}
	
\end{lem}

Proof. Since $3 \equiv 3 \pmod{4}$ if we apply the law of quadratic reciprocity we have two cases . If $n \equiv 1 \pmod{4}$ then $\left(\frac{3}{n}\right)=\left(\frac{n}{3}\right)$ and the result follows . If $n \equiv 3 \pmod{4}$ then $\left(\frac{3}{n}\right)=-\left(\frac{n}{3}\right)$ and the result follows .

\begin{lem}\label{le}
	Let $n$ be odd positive number , then	
	\begin{center}
	$\left(\frac{5}{n}\right)=
	\begin{cases}
	\hphantom{-}1&\text{if }n\equiv 1,4\pmod 5\\
	\hphantom{-}0&\text{if }n\equiv 0\pmod 5\\
	-1 & \text{if }n\equiv 2,3\pmod 5
	\end{cases}$
	\end{center}
	
\end{lem}

Proof. Since $5 \equiv 1 \pmod{4}$ if we apply the law of quadratic reciprocity we have $\left(\frac{5}{n}\right)=\left(\frac{n}{5}\right)$ and the result follows .

\begin{lem}\label{le}
	Let $n$ be odd positive number , then\\
	case 1. ($n\equiv 1 \pmod 4$)	
	\begin{center}
		$\left(\frac{-3}{n}\right)=
		\begin{cases}
		\hphantom{-}1&\text{if }n\equiv 1,11\pmod {12}\\
		\hphantom{-}0&\text{if }n\equiv 3,9\pmod {12}\\
		-1 & \text{if }n\equiv 5,7\pmod {12}
		\end{cases}$
	\end{center}
	case 2. ($n\equiv 3 \pmod 4$)
	\begin{center}
		$\left(\frac{-3}{n}\right)=
		\begin{cases}
		\hphantom{-}1&\text{if }n\equiv 5,7\pmod {12}\\
		\hphantom{-}0&\text{if }n\equiv 3,9\pmod {12}\\
		-1 & \text{if }n\equiv 1,11\pmod {12}
		\end{cases}$
	\end{center}
	
\end{lem}

Proof.$\left(\frac{-3}{n}\right)=\left(\frac{-1}{n}\right)\left(\frac{3}{n}\right)$ . Applying the law of quadratic reciprocity we have : if $n \equiv 1 \pmod {4}$ then $\left(\frac{3}{n}\right)=\left(\frac{n}{3}\right)$ . If $n \equiv 3 \pmod {4} $ then $\left(\frac{3}{n}\right)=-\left(\frac{n}{3}\right)$ . Applying the Chinese remainder theorem in both cases several times we get the result .

\begin{lem}\label{le}
	Let $n$ be odd positive number , then\\
	case 1. ($n\equiv 1 \pmod 4$)	
	\begin{center}
		$\left(\frac{7}{n}\right)=
		\begin{cases}
		\hphantom{-}1&\text{if }n\equiv 1,2,4\pmod {7}\\
		\hphantom{-}0&\text{if }n\equiv 0\pmod {7}\\
		-1 & \text{if }n\equiv 3,5,6\pmod {7}
		\end{cases}$
	\end{center}
	case 2. ($n\equiv 3 \pmod 4$)
	\begin{center}
		$\left(\frac{7}{n}\right)=
		\begin{cases}
		\hphantom{-}1&\text{if }n\equiv 3,5,6\pmod {7}\\
		\hphantom{-}0&\text{if }n\equiv 0\pmod {7}\\
		-1 & \text{if }n\equiv 1,2,4\pmod {7}
		\end{cases}$
	\end{center}
	
\end{lem}

Proof. Since $7 \equiv 3 \pmod{4}$ if we apply the law of quadratic reciprocity we have two cases . If $n \equiv 1 \pmod{4}$ then  $\left(\frac{7}{n}\right)=\left(\frac{n}{7}\right)$ and the result follows . If $n \equiv 3 \pmod{4}$ then $\left(\frac{7}{n}\right)=-\left(\frac{n}{7}\right)$ and the result follows .

\begin{lem}\label{le}
	Let $n$ be odd positive number , then\\
	case 1. ($n\equiv 1 \pmod 4$)	
	\begin{center}
		$\left(\frac{-6}{n}\right)=
		\begin{cases}
		\hphantom{-}1&\text{if }n\equiv 1,5,19,23\pmod {24}\\
		\hphantom{-}0&\text{if }n\equiv 3,9,15,21\pmod {24}\\
		-1 & \text{if }n\equiv 7,11,13,17\pmod {24}
		\end{cases}$
	\end{center}
	case 2. ($n\equiv 3 \pmod 4$)
	\begin{center}
		$\left(\frac{-6}{n}\right)=
		\begin{cases}
		\hphantom{-}1&\text{if }n\equiv 7,11,13,17\pmod {24}\\
		\hphantom{-}0&\text{if }n\equiv 3,9,15,21\pmod {24}\\
		-1 & \text{if }n\equiv 1,5,19,23\pmod {24}
		\end{cases}$
	\end{center}
	
\end{lem}

Proof.$\left(\frac{-6}{n}\right)=\left(\frac{-1}{n}\right)\left(\frac{2}{n}\right)\left(\frac{3}{n}\right)$ . Applying the law of quadratic reciprocity we have : if $n \equiv 1 \pmod {4}$ then $\left(\frac{3}{n}\right)=\left(\frac{n}{3}\right)$ . If $n \equiv 3 \pmod {4} $ then $\left(\frac{3}{n}\right)=-\left(\frac{n}{3}\right)$ . Applying the Chinese remainder theorem in both cases several times we get the result .

\begin{lem}\label{le}
	Let $n$ be odd positive number , then	
	\begin{center}
		$\left(\frac{10}{n}\right)=
		\begin{cases}
		\hphantom{-}1&\text{if }n\equiv 1,3,9,13,27,31,37,39\pmod {40}\\
		\hphantom{-}0&\text{if }n\equiv 5,15,25,35\pmod {40}\\
		-1 & \text{if }n\equiv 7,11,17,19,21,23,29,33\pmod {40}
		\end{cases}$
	\end{center}
	
\end{lem}

Proof.$\left(\frac{10}{n}\right)=\left(\frac{2}{n}\right)\left(\frac{5}{n}\right)$ . Applying the law of quadratic reciprocity we have :  $\left(\frac{5}{n}\right)=\left(\frac{n}{5}\right)$ . Applying the Chinese remainder theorem  several times we get the result .

\section{The Main Result}

\begin{thm}\label{le}
	Let $N=k\cdot 2^m-1$ such that $m>2$ , $3 \mid k$ , $0<k<2^m$ and 
	\begin{center}
		$\begin{cases}
		k \equiv 1 \pmod{10}~ with ~ m \equiv 2,3 \pmod{4} \\
		k \equiv 3 \pmod{10}~ with ~ m \equiv 0,3 \pmod{4} \\
		k \equiv 7 \pmod{10}~ with ~ m \equiv 1,2 \pmod{4} \\
		k \equiv 9 \pmod{10}~ with ~ m \equiv 0,1 \pmod{4}
		\end{cases}$ 
		\end{center}
		\begin{center}
		Let $b=3$ and $S_0(x)=V_k(b,1)$ , then \\
		$N$ is prime iff $N \mid S_{m-2}(x)$
	\end{center}
\end{thm}

Proof. Since $N \equiv 3 \pmod {4}$ and $b=3$ from Lemma 2.2 we know that $\left(\frac{2-b}{N}\right)=-1$ . Similarly , since $N \equiv 2 \pmod {5}$ or $N \equiv 3 \pmod {5}$ and $b=3$ from Lemma 2.5 we know that $\left(\frac{2+b}{N}\right)=-1$ . From Lemma 2.1 we know that $V_k(b,1)=x$ . Applying Theorem 2.1 in the case $c = 1$ we get the result.

\begin{thm}\label{le}
	Let $N=k\cdot 2^m-1$ such that $m>2$ , $3 \mid k$ , $0<k<2^m$ and 
	\begin{center}
	$\begin{cases}
	k \equiv 3 \pmod{42}~ with ~ m \equiv 0,2 \pmod{3} \\
	k \equiv 9 \pmod{42}~ with ~ m \equiv 0 \pmod{3} \\
	k \equiv 15 \pmod{42}~ with ~ m \equiv 1 \pmod{3} \\
	k \equiv 27 \pmod{42}~ with ~ m \equiv 1,2 \pmod{3} \\
	k \equiv 33 \pmod{42}~ with ~ m \equiv 0,1 \pmod{3} \\
	k \equiv 39 \pmod{42}~ with ~ m \equiv 2 \pmod{3} 
	\end{cases}$
	\end{center}
	 \begin{center}
	 Let $b=5$ and $S_0(x)=V_k(b,1)$ , then \\
	 $N$ is prime iff $N \mid S_{m-2}(x)$
	\end{center}
\end{thm}

Proof. Since $N \equiv 3 \pmod {4}$ and $N \equiv 11 \pmod {12}$ and $b=5$ from Lemma 2.6 we know that $\left(\frac{2-b}{N}\right)=-1$ . Similarly , since $N \equiv 3 \pmod {4}$ and $N \equiv 1 \pmod {7}$ or $N \equiv 2 \pmod {7}$ or $N \equiv 4 \pmod {7}$  and $b=5$ from Lemma 2.7 we know that $\left(\frac{2+b}{N}\right)=-1$ . From Lemma 2.1 we know that $V_k(b,1)=x$ . Applying Theorem 2.1 in the case $c = 1$ we get the result.

	\begin{thm}\label{le} 
		Let $N=k\cdot 2^m+1$ such that $m>2$ ,  $0<k<2^m$ and
		\begin{center}
			$\begin{cases}
			k \equiv 1 \pmod{42}~ with ~ m \equiv 2,4 \pmod{6} \\
			k \equiv 5 \pmod{42}~ with ~ m \equiv 3 \pmod{6} \\
			k \equiv 11 \pmod{42}~ with ~ m \equiv 3,5 \pmod{6} \\
			k \equiv 13 \pmod{42}~ with ~ m \equiv 4 \pmod{6} \\
			k \equiv 17 \pmod{42}~ with ~ m \equiv 5 \pmod{6} \\
			k \equiv 19 \pmod{42}~ with ~ m \equiv 0 \pmod{6} \\ 
			k \equiv 23 \pmod{42}~ with ~ m \equiv 1,3 \pmod{6} \\
			k \equiv 25 \pmod{42}~ with ~ m \equiv 0,2 \pmod{6} \\
			k \equiv 29 \pmod{42}~ with ~ m \equiv 1,5 \pmod{6} \\
			k \equiv 31 \pmod{42}~ with ~ m \equiv 2 \pmod{6} \\
			k \equiv 37 \pmod{42}~ with ~ m \equiv 0,4 \pmod{6} \\
			k \equiv 41 \pmod{42}~ with ~ m \equiv 1 \pmod{6}  
			\end{cases}$
		\end{center}
		\begin{center}
			Let $b=5$ and $S_0(x)=V_k(b,1)$ , then \\
			$N$ is prime iff $N \mid S_{m-2}(x)$
		\end{center}

	\end{thm}
	
	Proof. Since $N \equiv 1 \pmod {4}$ and $N \equiv 5 \pmod {12}$ and $b=5$ from Lemma 2.6 we know that $\left(\frac{2-b}{N}\right)=-1$ . Similarly , since $N \equiv 1 \pmod {4}$ and $N \equiv 3 \pmod {7}$ or $N \equiv 5 \pmod {7}$ or $N \equiv 6 \pmod {7}$  and $b=5$ from Lemma 2.7 we know that $\left(\frac{2+b}{N}\right)=-1$ . From Lemma 2.1 we know that $V_k(b,1)=x$ . Applying Theorem 2.1 in the case $c = 1$ we get the result.

\begin{thm}\label{le} 
	Let $N=k\cdot 2^m+1$ such that $m>2$ ,  $0<k<2^m$ and
	\begin{center}
		$\begin{cases} 
		k \equiv 1,7 \pmod{30} \text{ with } n \equiv 0 \pmod{4}  \\ 
		k \equiv 11,23 \pmod{30} \text{ with } n \equiv 1 \pmod{4}  \\ 
		k \equiv 13,19 \pmod{30} \text{ with } n \equiv 2 \pmod{4} \\ 
		k \equiv 17,29 \pmod{30} \text{ with } n \equiv 3 \pmod{4} 
		\end{cases}$
	\end{center}
	\begin{center}
		Let $b=8$ and $S_0(x)=V_k(b,1)$ , then \\
		$N$ is prime iff $N \mid S_{m-2}(x)$
	\end{center}

\end{thm}

Proof. Since $N \equiv 1 \pmod {4}$ and $N \equiv 17 \pmod {24}$ and $b=8$ from Lemma 2.8 we know that $\left(\frac{2-b}{N}\right)=-1$ . Similarly , since $N \equiv 17 \pmod {40}$ or $N \equiv 33 \pmod {40}$   and $b=8$ from Lemma 2.9 we know that $\left(\frac{2+b}{N}\right)=-1$ . From Lemma 2.1 we know that $V_k(b,1)=x$ . Applying Theorem 2.1 in the case $c = 1$ we get the result.

\makeatletter
\renewcommand{\@biblabel}[1]{[#1]\hfill}
\makeatother

\end{document}